\documentclass[12pt]{article}
\title{Monomial ideals and independence of $\mathrm{I}\Sigma_2$ \\ DRAFT}
\date{2015}
\author{Florian Pelupessy\footnote{Partially supported by the Japan Society for the Promotion of Science (KAKENHI 23340020) and a JSPS postdoctoral fellowship for foreign researchers.
This is the pre-peer reviewed version of an article which is in press at Mathematical Logic Quarterly} }
\usepackage[english]{babel}
\usepackage{amsmath, amssymb}
\usepackage{url}
\usepackage{xcolor}
\usepackage{hyperref}
\usepackage{tikz}
\usepackage{multicol}
\usepackage{parskip}
\usepackage{makeidx}
  \makeindex

\usepackage[initials]{amsrefs}

\BibSpec{article}{%
  +{}{\PrintAuthors}				{author} %
  +{,}{ \textit}					{title}
  +{,}{ }							{journal}
  +{}{ \textbf}					{volume}
  +{,}{ }							{note}
  +{}{ \parenthesize}				{date}
  +{:}{ }							{pages}
  +{.}{}							{transition}
  +{}{ }							{review}
  +{\\}{ \url }							{url}
}

\BibSpec{thesis}{%
  +{}{\PrintAuthors}				{author} %
  +{,}{ \textit}					{title}
  +{,}{ }							{journal}
  +{}{ \textbf}					{volume}
  +{,}{ }							{note}
  +{,}{ }							{date}
  +{:}{ }							{pages}
  +{.}{}							{transition}
  +{}{ }							{review}
  +{\\}{ \url }							{url}
}

\usepackage[T1]{fontenc}
\usepackage{libertine}

\newtheorem{theorem}{Theorem}
\newtheorem{lemma}{Lemma}
\newtheorem{definition}{Definition}

\newtheorem{corollary}[theorem]{Corollary}

\begin{document}
\maketitle
\begin{abstract}
We show that a miniaturised version of Maclagan's theorem on monomial ideals is equivalent to $\mathrm{1}{-}\mathrm{Con}(\mathrm{I}\Sigma_2)$ and classify a phase transition threshold for this theorem. This work highlights the combinatorial nature of Maclagan's theorem.
\end{abstract}

Monomial ideals play an important role in commutative algebra and algebraic combinatorics. Maclagan's theorem has several applications in computer algebra (See, for example, \cite{maclagan}), so the logical and combinatorial issues which surround it are of particular interest. 

We first determine an upper bound for the miniaturised version of Maclagan's theorem based on a very short proof of Maclagan's theorem using Friedman's adjacent Ramsey theorem. We apply known upper bounds for Friedman's finite adjacent Ramsey theorem for this part.

We provide lower bounds for the miniaturised Maclagan theorem. This proof of independence of $\mathrm{I}\Sigma_2$ is a cleaned up version of the proof in~\cite{pelupessyweiermann}. We finish by determining a sharp phase transition threshold for the miniaturised Maclagan's theorem. 

These results complement the study by Aschenbrenner and Pong~\cite{aschpong} and the determination of lower bounds in~\cite{pelupessyweiermann}. Furthermore this paper fits into the general programme started by Andreas Weiermann (see~\cite{weiermannweb}) to determine the phase transition thresholds for a variety of theorems. The work by Aschenbrenner and Pong demonstrates that the infinite version of Maclagan's theorem's strength is at the level of the well ordering of $\omega^{\omega^\omega}$. The determination of the equivalence to $\mathrm{1}{-}\mathrm{Con}(\mathrm{I}\Sigma_2)$ shows the logical strength of the miniaturised version, which is at the same level as the Paris--Harrinton theorem for triples, miniaturised Higman's lemma and the finite adjacent Ramsey theorem for pairs. The phase transition for Maclagan's theorem is the same is that of the latter two of those theorems. 

\section{Preliminaries}
In what follows we take arbitrary field $K$ and examine ideals in the polynomial ring $K[X_d, \dots, X_0, Y]$. A monomial \index{monomial} is a polynomial of the form $m=X_d^{i_d} \dots X_0^{i_0}Y^j$. We use $\bar{m}$ to denote the corresponding $d+2$-tuple $(i_d, \dots, i_0,j)$. For $d+2$-tuples $x \leq y$ denotes the coordinatewise ordering. A monomial ideal is an ideal that is generated by monomials. We denote an ideal that is generated by a set $G$ of generators with $\langle G \rangle$ \index{$\langle G \rangle$}. Notice that if $m$ is a monomial then $m \in \langle G \rangle$ if and only if there is an $m' \in G$ such that $\bar{m}' \leq \bar{m}$. 
\begin{theorem}[Maclagan]\label{thm:Maclagan}
Every infinite sequence of monomial ideals contains an ideal that is a subset of an ideal that occurs earlier in the sequence.
\end{theorem}
We prove this with the infinite adjacent Ramsey theorem for pairs with codimension $d+2$:
\begin{theorem}[Adjacent Ramsey]
For every $C\colon \mathbb{N}^2 \rightarrow \mathbb{N}^r$ there exist $x_1 < x_2 < x_3$ such that $C(x_1,x_2) \leq C(x_2,x_3)$.
\end{theorem}
\emph{Proof:} From~\cite{friedman2010}: Use Ramsey's theorem for triples with $r+1$ colours on the colouring:
\[
D(x,y,z)=
\left\{ 
\begin{array}{ll}
0 & \textrm{if $C(x,y) \leq C(y,z)$} \\
i & \textrm{otherwise},
\end{array}
\right.
\]
where $i$ is the least such that $(C(x,y))_i > (C(y,z))_i$.
\begin{flushright} $\Box$ \end{flushright}
\emph{Proof of \ref{thm:Maclagan}:} Suppose, for a contradiction, that we have an infinite sequence $I_0, I_1, I_2, \dots$ of monomial ideals with $I_i \not\supseteq I_j$ for all $i<j$, where $I_i=\langle \{ m_i^1 , \dots , m_i^{n_i} \} \rangle$. Define:
\[
\mathrm{P}(i,j)=\mathrm{the \ smallest \ } k \mathrm{\ such \ that \ } \forall l\leq n_i.\bar{m}_i^l \not\leq \bar{m}_j^k.
\]
Take:
\[
C(i,j)=\bar{m}_j^{\mathrm{P}(i,j)}. 
\]
By adjacent Ramsey there exist $a<b<c$ such that $C(a,b) \leq C(b,c)$, contradiction!
\begin{flushright} $\Box$ \end{flushright}
\begin{definition} ~{}
\begin{enumerate}
\item The degree of a monomial is the total degree:
\[ \mathrm{deg}(X_d^{i_d} \dots X_0^{i_0}Y^j)=i_d+ \dots + i_0+j.\]
\item The degree of a finite set $G$ of monomials is the maximum of the degrees of the elements of that set: $\mathrm{deg}(G)=\max \{ \mathrm{deg}(m): m \in G \}$.
\item The degree of a monomial ideal $I$ is the smallest degree that is needed to be able to generate it with monomials: $\mathrm{deg}(I)= \min \{ \mathrm{deg}(G): I=\langle G \rangle \}$.
\end{enumerate}
\end{definition}

\begin{theorem}[$\mathrm{MM}$, miniaturised Maclagan] \label{thm:MM}
For every $d,l$ there exists an $M$ such that for every sequence $I_0, \dots, I_M$ of monomial ideals in $K[X_d, \dots, X_0, Y]$,   with $\mathrm{deg}(I_i) \leq l+i$ for all $i \leq M$, there exist $i < j \leq M$ with $I_i \supseteq I_j$.
\end{theorem}
\emph{Proof:} Repeat the proof of Theorem~\ref{thm:Maclagan} with the following version of the finite adjacent Ramsey theorem:
\begin{theorem}
For every $d,l$ there exists $R$ such that for every  
\[C\colon \{l, \dots ,R\}^2 \rightarrow \mathbb{N}^{d+2}
\]
with $\max C(x,y) \leq \max \{ x,y\}$ for all $(x,y) \in \{l , \dots , R\}^2$ there exist $a<b<c$ with $C(a,b) \leq C(b,c)$. 
\end{theorem}
\emph{Proof:} Shown from the finite adjacent Ramsey theorem from~\cite{friedmanpelupessy} with increased codimension.
\begin{flushright} $\Box$ \end{flushright}
\begin{corollary}
$\mathrm{EFA} \vdash \mathrm{1}{-}\mathrm{Con}(\mathrm{I}\Sigma_2) \rightarrow \mathrm{MM}$
\end{corollary}
\emph{Proof:} Use the upper bounds for the finite adjacent Ramsey theorem from~\cite{friedmanpelupessy}.
\begin{flushright} $\Box$ \end{flushright}
We introduce a parameter $f\colon \mathbb{N} \rightarrow \mathbb{N}$:
\begin{theorem}[$\mathrm{MM}_f$]
For every $l$ there exists an $M$ such that for every sequence $I_0, \dots, I_M$ of monomial ideals in $K[X_d, \dots, X_0,Y]$,   with $\mathrm{deg}(I_i) \leq l+f(i)$ for all $i \leq M$, there exist $i < j \leq M$ with $I_i \supseteq I_j$.
\end{theorem}
\emph{Proof:} Compactness on the infinite Maclagan's theorem.
\begin{flushright} $\Box$ \end{flushright}
We use $M_d^f(l)$ \index{$M_d^f$} to denote the least $M$ from $\mathrm{MM}_f$. We assume basic knowledge of ordinals $\leq \omega^{\omega}$, their cantor normal forms and the canonical fundamental sequences. Furthermore we use the fast growing (or Wainer) hierarchy:
\begin{eqnarray*}
F_0(i) & = & i+1, \\
F_{\alpha+1}(i) & = & F_\alpha^i(i), \\
F_\gamma (i) & = & F_{\gamma[i]} (i). 
\end{eqnarray*}
and the fact that any computable function is provably total in $\mathrm{I}\Sigma_2$ if and only if it is primitive recursive in $F_\alpha$ for some $\alpha<\omega^{\omega}$. More information on ordinals and fast growing hierarchies can be found in \cite{araiord}, \cite{rose}, information on the connections between hierarchies and provability in \cite{takeuti} or \cite{buchholz, arai, pohlers, buss}.

All parameter functions are assumed to be nondecreasing. For every unbounded $f\colon \mathbb{N} \rightarrow \mathbb{N}$ the inverse is:
\[
f^{-1} (i) = \max \{ j: f(j) \leq i\}\cup\{0\}.
\]
$\log$ is the inverse of $i \mapsto 2^i$, $\log \log$ is the inverse of $i \mapsto 2^{2^i}$, $\sqrt[c]{\log}$ is the inverse of $i \mapsto 2^{(i^c)}$ and $i \mapsto \frac{i}{c}$ is the inverse of $i \mapsto i \cdot c$. 
The remainder of this article is dedicated to proving the following:
\begin{theorem}
Define $f_\alpha(i)=\sqrt[F^{-1}_\alpha (i)]{\log i}$. 
\begin{enumerate}
  \item $\mathrm{I}\Sigma_2 \nvdash \mathrm{MM}_{\mathrm{id}}$.
  \item $\mathrm{I}\Sigma_2 \nvdash \mathrm{MM}_{\sqrt[c]{\log}}$.
  \item $\mathrm{I}\Sigma_2 \vdash \mathrm{MM}_{\log \log}$.
  \item $\mathrm{I}\Sigma_2 \vdash \mathrm{MM}_{f_\alpha}$ if and only if $\alpha<\omega^{\omega}$.
\end{enumerate}
\end{theorem}
Item 1 is Theorem~\ref{thm:MM:id}, item 2 and the unprovability part of item 4 are Theorems~\ref{thm:MM:rootlog} and~\ref{thm:MM:slowbutunprovable} and item 3 and the provability part of item 4 are shown with Lemma~\ref{lemma:MM:slowandprovable}. 

\section{Lower bounds for the identity function}
We will first examine the independence of $\mathrm{MM}_\mathrm{id}$. Compared with the proof in \cite{pelupessyweiermann} we have removed the inconvenient step of constructing special descending sequences of ordinals and encoded those directly in the monomial ideals. Furthermore intermediate sequences have been removed, greatly simplifying the proof (especially the bookkeeping of the degrees of the ideals in sequences).
\begin{theorem}\label{thm:MM:id}
$\mathrm{I}\Sigma_2 \nvdash \mathrm{MM}_\mathrm{id}$
\end{theorem}
\emph{Proof:}
We call a sequence $I_0, \dots, I_R$ with $\mathrm{deg}(I_i) \leq l+i$ for all $i \leq R$ bad \index{bad!sequence of monomial ideals} if there do not exist $i < j \leq R$ with $I_i \supseteq I_j$. We call a sequence of sets of generators bad if the sequence of ideals generated by those sets is bad.  Such a sequence shows that $R<M_d^\mathrm{id} (l)$.

We will associate with each ordinal $\alpha<\omega^{d+1}$ some monomials and number $h_\alpha$. We construct sequences of sets of generators consisting of monomials associated with ordinals $\leq \alpha$ that show $F_\alpha(l) < M_d^\mathrm{id} (l+h_\alpha)$. We leave out the brackets in the definitions of the sets of generators and identify ideals with their sets of generators. 

Given:
\[
\alpha=\omega^d \cdot n_d + \dots + \omega^0 \cdot n_0,
\]
we associate with $\alpha$ the set of monomials of the form:
\[
X_d^{n'_d} \cdots X_0^{n'_0}Y^m,
\]
where $n'_i \leq 2 \cdot n_i+1$, $m \in \mathbb{N}$ and the number $h_\alpha=2\cdot n_d + \dots + 2 \cdot n_0 + d+1$. Because the existence of such a sequence implies:
\[
F_{\omega^d}(l) < M_d^{\mathrm{id}}(l+3+d+1),
\]
for all $d$ this construction suffices to prove Theorem~\ref{thm:MM:id}. 

The bad sequences are defined by recursion on $\alpha$:
\begin{itemize}
\item
For $\alpha=0$ we take the following sequence:
\begin{eqnarray*}
\mathrm{Seq}(\alpha, l)_0 & = & X_0, \\
\mathrm{Seq}(\alpha, l)_{1+i} & = & Y^{l+1-i} \textrm{ \ \ for $0 \leq i \leq l+1$}.
\end{eqnarray*}
\item
For $\alpha+1$ we start the construction with, for $0 \leq i \leq l$:
\[
\mathrm{Seq}(\alpha+1, l)_i  =  X_d^{2 \cdot n_d+1} \cdots X_1^{2 \cdot n_1+1}X_0^{2 \cdot n_0+3}Y^{l-i}.
\]
Continuing with, for $F_\alpha^0(l) + \dots  +F_\alpha^j(l) < i \leq F_\alpha^0(l)+ \dots +F_\alpha^{j+1} (i)$ and $0\leq j \leq l$:
\[
\mathrm{Seq}(\alpha+1, l)_{i}  =  X_d^{2 \cdot n_d+1} \cdots X_1^{2\cdot n_1+1}X_0^{2 \cdot n_0+2}Y^{l-j}, b_i \cdot Y^{l+1},
\]
where $b_i=\mathrm{Seq}(\alpha, F_\alpha^j(l))_{i-F_\alpha^0(l) - \dots  -F_\alpha^j(l)}$ and $b_im$ denotes the set consisting of the elements from $b_i$ multiplied by monomial $m$. 
\item
For limit $\alpha=\omega^d \cdot n_d + \dots + \omega^j \cdot (n_j+1)+ \dots +\omega^0 \cdot n_0$ (where $j>0$ and $n_{j-1}, \dots, n_0=0$) we take the sequence defined as follows, for $0 \leq i \leq l$:
\[
\mathrm{Seq}(\alpha, l)_i  =  X_d^{2\cdot n_d+1} \cdots X_{j+1}^{2 \cdot n_{j+1} +1}X_0^{2 \cdot n_j +3}Y^{l-i}. 
\]
For $0 \leq i \leq 2l$:
\[
\mathrm{Seq}(\alpha, l)_{l+i+1} = X_d^{2\cdot n_d+1} \cdots X_{j+1}^{2 \cdot n_{j+1} +1}X_0^{2 \cdot n_j +2}Y^{2l-i}.
\]
For $0 < i \leq F_{\alpha[l]}(l)$:
\[
\mathrm{Seq}(\alpha, l)_{3l+i+1}  =  \mathrm{Seq}(\alpha[l], l)_{i}.
\]
\end{itemize}
\emph{Claim:} The sequences $\mathrm{Seq}(\alpha, l)$ show $F_\alpha(l) < M_d^\mathrm{id} (l+h_\alpha)$. 

\emph{Proof:} According to the definition the sequences $\mathrm{Seq}(\alpha, l)$ are long enough. First we show, using induction on $\alpha$, that the degrees are bounded linearly:
\begin{itemize}
\item $\alpha=0$: $\mathrm{deg}(X_0)=1$ and, if $0<i\leq l+2$ then we have:
\[
 \mathrm{deg}(\mathrm{Seq}(\alpha, l)_{i})=l+2-i \leq l\leq l+h_0+i.
\]
\item $\alpha+1$: if $0 \leq i \leq l$ then:
\[
\mathrm{deg}(\mathrm{Seq}(\alpha, l)_i)=2\cdot n_d+ \dots +2\cdot n_0+d+3+l-i \leq l+h_{\alpha+1}+i.
\]
If $F_\alpha^0(l) + \dots  +F_\alpha^j(l) < i \leq F_\alpha^0(l)+ \dots +F_\alpha^{j+1} (i)$ and $0\leq j \leq l$ then:
\[
 \mathrm{deg}(\mathrm{Seq}(\alpha, l)_{i})\leq \max \{h_{\alpha+1}, \mathrm{deg}(b_i)+l+1\} \leq l+h_{\alpha+1}+i,
\]
where the second inequality is obtained from the definition of $b_i$ and the induction hypothesis.
\item Limit $\alpha$: If $0 \leq i \leq l$ then:
\[
\mathrm{deg}(\mathrm{Seq}(\alpha, l)_i)\leq l+ h_\alpha -i \leq l+h_\alpha+i.
\]
If $l<i\leq 3l+1$ then:
\[
\mathrm{deg}(\mathrm{Seq}(\alpha, l)_i)\leq l+h_\alpha + (l+1-i) \leq l+h_\alpha+i.
\]
If $3l+1< i \leq 3l+1 + F_{\alpha [l]}(l)$ then:
\[
\mathrm{deg}(\mathrm{Seq}(\alpha, l)_i)= \mathrm{deg}(\mathrm{Seq}(\alpha[l], l)_{i-3l-1}) \leq l+h_{\alpha [l]}+ i-3l-1,
\]
where the latter inequality results from the induction hypothesis. Notice that $h_{\alpha [l]}=h_\alpha+2l$, hence:
\[
\mathrm{deg}(\mathrm{Seq}(\alpha, l)_i) \leq l+h_\alpha+2l+i-3l-1 \leq h_\alpha+i-1 \leq l+h_\alpha+i.
\]
\end{itemize} 
We still need to prove that these sequences are bad. We will use the fact that a monomial is an element of a monomial ideal if and only if one of the generators divides that monomial and that in the construction the generators in each sequence consist of monomials which are associated with ordinals $\leq \alpha$ exclusively.
\begin{itemize}
\item $\alpha=0$: Notice that $X_0$ does not divide $Y^a$ and if $i<j\leq l+1$ then $Y^{l+1-i}$ does not divide $Y^{l+1-j}$.
\item $\alpha+1$: The generators $X_d^{2 \cdot n_d+1} \cdots X_1^{2 \cdot n_1+1}X_0^{2 \cdot n_0+1+a}Y^{b}$ ($a=1,2$) do not divide any monomials that are associated with $\beta \leq \alpha$. Indeed, if such a generator divided such a monomial we would have $2\cdot n_0+1+a \leq 2\cdot n_0+1$. Hence if 
\[
F_\alpha^0(l)+ \dots + F_\alpha^{j}(l)<i_0 <i_1\leq F_\alpha^0(l)+ \dots + F_\alpha^{j+1}(l),
\]
then any generator of $\mathrm{Seq}(\alpha+1,l)_{i_1}$ that is associated with ordinals $\leq \alpha$ is not divided by $X_d^{2 \cdot n_d+1} \cdots X_1^{2 \cdot n_1+1}X_0^{2 \cdot n_0+2}Y^{l-j}$. By induction hypothesis those generators (from $b_{i_1}Y^{l+1}$) can also not be divided by the other elements of $\mathrm{Seq}(\alpha+1,l)_{i_0}$. 

If
\[
i_0\leq F_\alpha^0(l)+ \dots + F_\alpha^{j}(l) < i_1,
\]
then the generator $X_d^{2 \cdot n_d+1} \cdots X_1^{2 \cdot n_1+1}X_0^{2 \cdot n_0+1+a}Y^{b}$ is not divided by any element of $\mathrm{Seq}(\alpha+1,l)_{i_0}$.
\item Limit $\alpha$: Again, if $i_0< 3l+1<i_1$ then the generators in $\mathrm{Seq}(\alpha,l)_{i_1}$ are associated with ordinals $<\alpha$, hence cannot be divided by any generator from $\mathrm{Seq}(\alpha,l)_{i_0}$. If $3l+1<i_0<i_1$ then the induction hypothesis delivers the same fact as does the definition of $\mathrm{Seq}(\alpha,l)$ when $i_0<i_1\leq l+1$ or $l+1< i_0<i_1 \leq 3l+1$.
\end{itemize}
This ends the proof of the claim thus finishing the proof Theorem~\ref{thm:MM:id}.
\begin{flushright} $\Box$ \end{flushright}

\begin{corollary}
$\mathrm{EFA} \vdash \mathrm{1}{-}\mathrm{Con}(\mathrm{I}\Sigma_2) \leftrightarrow \mathrm{MM}_{\mathrm{id}}$
\end{corollary}

\section{Lower bounds for other parameter values}
In this section we modify bad sequences for the identity into bad sequences for lower parameter values $f$, showing that $M^f$ again is unbounded in the multiply recursive functions. The first step of the modification is to, given a sequence:
\[
I_0, \dots ,I_M,
\]
define the new sequence: 
\[
I_{f(0)}, \dots , I_{f(M)}.
\]
This new sequence will have identical elements, but we correct this by modifying the ideals using a constant number of extra variables. To be able to do this we will need to estimate the number $c'$ sufficiently large such that
\[
\# \{ i : f(i)=f(j) \} \leq M_{c'}^0 (j). 
\]
For this reason we start with studying $M^0$.
\begin{lemma}\label{lemma:MM:zero:exp}
$M_0^0(2j+2) \geq 2^j$.
\end{lemma}
\emph{Proof:} We construct sequences that show this using recursion on $j$. 

For $j=0$ we take sequence $m_0=Y^2$, $m_1=Y$. 

Given sequence $a_0, \dots , a_{2^j}$ for $j$, take for $j+1$ the sequence defined by:
\[
m_i = \left\{ 
\begin{array}{ll}
X_0^{j+1}Y^1, a_iY^2 & \textrm{if $0 \leq i \leq 2^j$}, \\
X_0^{j+1}Y^0, a_{i-2^j-1}Y^2 & \textrm{if $2^j+1 \leq i \leq 2^{j+1}$}. \\
\end{array}
\right.
\]
To show that this is a desired sequence notice first that $\mathrm{deg}(a_i)\leq 2j-1$. Hence $\mathrm{deg}(m_i) \leq 2j+2$ and $\langle m_i \rangle \not \supseteq \langle m_{i'} \rangle$ for $i<i'$.
\begin{flushright} $\Box$ \end{flushright}
\begin{lemma}\label{lemma:MM:zero:polying}
$M_{d+2}^0 (j+c) \geq M_d^0 (j)^{c+1}$.
\end{lemma}
\emph{Proof:} We construct sequences which show this. Taking a bad sequence $a_0, \dots a_M$, the elements of the new sequence will look like:
\[
X_{d+2}^cX_{d+1}^0a_{i_0}, \dots, X_{d+2}^{c-j}X_{d+1}^ja_{i_j}, \dots , X_{d+2}^0X_{d+1}^ca_{i_c}.
\]
The main idea is that the generators of the ideals in the new sequence get separated into `tracks', where due to the part $X_{d+2}^{c-j}X_{d+1}^j$ the generators from different tracks cannot divide each other. Hence if we change in this set of generators $a_{i_j}$ into $a_{i_j+1}$ the ideal generated by this new set is not a subset of the original ideal. Using this we construct the new sequences by recursion on $c$. 

For $c=0$ we take $m_i=a_i$. \\
Given sequence $b_0, \dots , b_N$ for $c$ and $0 \leq i \leq M \cdot N$ we take:
\[
m_i = X_{d+2}a_{\frac{i}{N}}, X_{d+1}b_{i \% N}. 
\]
To show this is a desired sequence notice first that $\mathrm{deg}(m_i)\leq \max \{a_i,b_i\} +1 \leq j+c+1$. To show that this is a bad sequence examine the following cases for $i<j$:
\begin{itemize}
  \item $\frac{i}{N}=\frac{j}{N}$ and $i \% N < j \% N$: Because $X_{d+1}b_{i \% N}$ does not contain $X_{d+2}$, $\langle m_i \rangle \not\supseteq \langle m_j \rangle$ is inherited from the $b_i$'s.
  \item $\frac{i}{N}< \frac{j}{N}$: Because $X_{d+2}a_{\frac{i}{N}}$ does not contain $X_{d+1}$, $\langle m_i \rangle \not\supseteq \langle m_j \rangle$ is inherited from the $a_i$'s.  
\end{itemize}
\begin{flushright} $\Box$ \end{flushright}
\begin{lemma}\label{MM:zero:polypower}
$M_{2c}^0(2^{c+1}(j+1)) \geq 2^{j^{c+1}}$.
\end{lemma}
\emph{Proof:} Combine the previous two lemmas, starting with the first for $c=0$, for the induction step we use the latter.
\begin{flushright} $\Box$ \end{flushright}
With the constructions so far it is not possible to obtain double-exponential lengths of such sequences, an attempt to do so would require using a non-constant number of variables. We will later see that double exponential lengths are not possible using any construction due to the upper bounds on $M^0$. We take this `highest possible' estimate of $M^0$ to prove unprovability for the following `low' parameter.
\begin{theorem}\label{thm:MM:rootlog}
If $f_c(i)=\sqrt[c]{\log (i)}$ then: 
$\mathrm{I}\Sigma_2 \nvdash \mathrm{MM}_{f_c}$.
\end{theorem}
\emph{Proof:} We use Lemma~\ref{MM:zero:polypower} to convert bad sequences for identity into bad sequences for $f_c$. Together with Theorem~\ref{thm:MM:id} the following is sufficient to prove this theorem:   
\[
M_{d+2c+3}^{f_c}(l) \geq M_d^{\mathrm{id}}(l)=:M,
\]
for $l \geq 2^{(c+4)^2}+1$. Our building blocks are bad sequences $a_0, \dots, a_M$ from the proof of Theorem~\ref{thm:MM:id} ($d+2$ variables) and $b(i)_0, \dots b(i)_{2^{(i+1)^c}}$ from Lemma~\ref{MM:zero:polypower} ($2c+2$ new variables and $i > 2^{(c+1)^2}$). The new sequence is:
\[
m_i = \left\{ 
\begin{array}{ll}
X_{d+2c+3}^{2^{(c+1)^2}-i+1} & \textrm{if $i \leq 2^{(c+1)^2}$}\\
a_{f_c(i)}, b(f_c(i))_{i-2^{f_c(i)^c}}& \textrm{otherwise}. \\ 
\end{array}
\right.
\]
First note that $\mathrm{deg}(m_i)\leq \max \{ 2^{(c+1)^2}, l+f_c(i), l+f_c(i) \} \leq l+f_c(i)$. It remains to show that this sequence is bad, for $i<j$ we have the following cases:
\begin{itemize}
  \item $i<j\leq 2^{(c+1)^2}$: $\langle m_i \rangle \not\supseteq \langle m_j \rangle$ follows directly from the definition.
  \item $i \leq 2^{(c+1)^2}< j$: The generators of $m_i$ contain $X_{d+2c+3}$ whilst those of $m_j$ do not, hence $\langle m_i \rangle \not\supseteq \langle m_j \rangle$.
  \item $2^{(c+1)^2}< i<j$: The generators of the $a$'s contain different variables from the generators of the $b$'s, hence $\langle m_i \rangle \not\supseteq \langle m_j \rangle$ is inherited from the $a$'s whenever $f_c(i) < f_c(i)$ and inherited from the $b$'s if $f_c(i)=f_c(j)$.
\end{itemize}
\begin{flushright} $\Box$ \end{flushright}
\begin{theorem}\label{thm:MM:slowbutunprovable}
If $f(i)=\sqrt[F_{\omega^\omega}^{-1}(i)]{\log (i)}$ then: 
$\mathrm{I}\Sigma_2 \nvdash \mathrm{MM}_f$.
\end{theorem}
\emph{Proof:} We show:
\[
M^f_{3l+3}(2^{(l+4)^2}+2l+4) \geq F(l).
\]
Assume for a contradiction that $M=M^f_{3l+3}(2^{(l+4)^2}+2l+4) < F(l)$. For $i \leq M$ we know that $F^{-1}(i) \leq l$, in other formulas that $\sqrt[F^{-1}(i)]{\log (i)} \geq \sqrt[l]{\log (i)}=f_l(i)$, using notation from the previous theorem. The estimates from Theorems~\ref{thm:MM:id} and ~\ref{thm:MM:rootlog} deliver:
\begin{eqnarray*}
M & \geq & M^{f_l}_{3l+3}(2^{(l+4)^2}+2l+4) \\
  & \geq & M^{\mathrm{id}}_l(2l+4) \\
  & \geq & F_{\omega^l} (l)=F(l)
\end{eqnarray*} 
which contradicts our assumption.
\begin{flushright} $\Box$ \end{flushright}
\section{Upper bounds}
For the upper bounds we use a simple counting of monomials, notice first that:
\begin{lemma} \label{MM:c}
$\mathrm{EFA} \vdash \mathrm{MM}_c$ for constant function $c$.
\end{lemma}
\emph{Proof:} The number of ideals of degree less than $l+c$ is bounded by the number of sets of monomials of degree less than $l+c$. This number has (rough) upper bound $2^{(l+c+1)^{d+2}}$. Hence, by pigeon hole principle, $M_d^c(l) \leq 2^{(l+c+1)^{d+2}}+1$.
\begin{flushright} $\Box$ \end{flushright}
%Using essentially the same argument we obtain the following result:
\begin{lemma}\label{lemma:MM:slowandprovable}
If $B$ is an increasing multiply recursive function and 
\[
f(i)=\sqrt[B^{-1}(i)]{\log (i)}
\]
then:
\[
\mathrm{I}\Sigma_2 \vdash \mathrm{MM}_f.
\] 
\end{lemma}
\emph{Proof:} Assume without loss of generality that $B(l)>2^{l+2}$. We show for $l>d+2$:
\[
M_d^f(l) < 2^{B(l)^{d+2}}.
\]
Take $R=2^{B(l)^{d+2}}$ and any sequence of monomial ideals $I_0, \dots , I_R$ with $\mathrm{deg}(I_i) \leq l+f(i)$. In this case we have:
\[
\mathrm{deg}(I_i) \leq l+\sqrt[B^{-1}(R)]{\log (R)} \leq l+B(l)^{\frac{d+2}{l}}\leq l+B(l).
\]
So by the upper bounds from the proof of lemma \ref{MM:c}:
\[
M_d^f(l) \leq M_d^{B(l)}(l) \leq 2^{(l+B(l)+1)^{d+1}}+1 < 2^{2^{d+1}B(l)^{d+1}}+1 <2^{B(l)^{d+2}}.
\]
\begin{flushright} $\Box$ \end{flushright}

\end{document}